# Noncommmutative Gelfand Duality and Applications I: The existence of invariant subspaces

Mukul S. Patel
<patel@math.uga.edu>
Department of Mathematics, University of Georgia
Athens, GA 30602, USA

February 13, 2006

**Abstract**

Gelfand duality between unital commutative $C^*$-algebras and Compact Hausdorff spaces is extended to all unital $C^*$-algebras, where the dual objects are what we call compact Hausdorff *quantum* spaces. We apply this result to obtain a characterization of unitary groups of $C^*$-algebras, and, for arbitrary bounded Hilbert space operators, (i) A spectral theorem *cum* continuous functional calculus, and (ii) A proof of the general Invariant Subspace Theorem. Also described is a nonabelian generalization of Pontryagin duality of abelian locally compact groups.

## 1   INTRODUCTION

(*All algebras and operators considered are over $\mathbb{C}$, the field of complex numbers, and so are all functions. Also, unless otherwise stated, all algebras considered are unital.*)

### 1.1   Description of the article

Connes attaches $C^*$-algebras to various geometrical objects arising in wide range of mathematical contexts [8], so that the geometry is studied with largely algebraic methods. Conversely, we assign a natural 'quantum space' to any given $C^*$-algebra, and use the former to study the latter. Of course, for commutative algebras, the Gelfand-Naimark theorem does the job:

**Theorem 1 (Gelfand-Naimark)**  *A commutative unital $C^*$-algebra $A$ is naturally isomorphic to $C(P(A))$, the algebra of complex-valued continuous functions on $P(A)$, the space of pure states of $A$.*

Note that $P(A)$ is a compact Hausdorff space, and the theorem sets up a functorial equivalence between the category of unital commutative $C^*$-algebras on one hand and the category of compact Hausdorff spaces on the other hand.

Since the appearance of Theorem 1 [14], there have been several noncommutative generalizations in various directions [1, 2, 3, 4, 24, 6, 11, 12, 20] with varying degree of success. Two works closest in spirit and content to the present article are elaborated in [3, 4, 24]. We will compare these with our approach after we explain our approach. For now, we note that applications of above generalizations have not been as extensive as those of the commutative Gelfand-Naimark theorem. Our generalization completely codifies the algebra structure into a topological object, and seems to be more natural. Ultimately, though, it is the range of applications that seems to justify our generalization. The latter (Theorem 3) is implemented by identifying the *natural* noncommutative topological analog of compact Hausdorff space, *i.e.,* what we call a compact Hausdorff *quantum space* (Definiton 7). Before we describe this object, we emphasize the following:

*Our notion of a quantum space is quite different than many hitherto considered.*

While it is a truism that the noncommutativity of a $C^*$-algebra $A$ is encoded in the non-Hausdorffness of its spectrum $Sp(A)$, we are not taking the latter as the topological representative of the algebra. Instead, our point of departure (Proposition 1) is the following circumstance: *A unital $C^*$-algebra $A$ is commutative if and only if all its irreducible Gelfand-Naimark-Segal representations are pair-wise inequivalent.* Thus, the noncommutativity of an algebra is completely captured by the equivalence relation given by equivalence of irreducible GNS representations. Equivalently, this gives an equivalence relation on $P(A)$, the space of pure states, because a state is pure if and only if the corresponding GNS representation is irreducible. We denote this equivalence relation on $P(A)$ by $R(A)$. In Section 2 are identified the primary obstacles to naive attempts at extending Gelfand-Naimark duality to noncommutative algebras.

Fortunately enough, these obstacles can be overcome by using the algebra $A$ to endow $P(A)$ with a natural *quantum space* structure (Definiton 7), in which the equivalence relation $R(A)$ is a crucial structural element.

Towards this end, we introduce a notion of *quantum sets* in Section 3. Again, our quantum notion is much different, direct, and concrete compared to other proposed notions of quantum sets, for example, as in [26]. Next is defined a natural noncommutative product of functions on a quantum set. We refer to this product as *q-product*. Next, in Section 4, we define the notion of *quantum topology,* or briefly *q-topology,* on a quantum set. Likewise, a quantum set with a q-topology will be abreviated as *q-space,* and functions continuous with respect to a q-topology will be called *q-continuous.* Then we show that the q-product of q-continuous functions on a quantum space is again a q-continuous function. Also demonstrated is the basic fact that given a compact Hausdorff q-space $X$, the set $A(X)$ of q-continuous functions on $X$, equipped with the q-product, is a unital $C^*$-algebra under the sup-norm.

We return to the $C^*$-algbera $A$ in Section 5, where the structure of $A$ is used to define a natural q-space structure on $P(A)$ making it a compact Hausdorff q-space. Instead of the spectrum $Sp(A)$, we take the q-space $P(A)$ as the topological 'dual' of $A$. Then, $Sp(A)$ of $A$ is the *quotient* of $P(A)$ under the equivalence relation $R(A)$. As above, the noncommutative product on space $A(P(A))$ of q-continuous functions on $P(A)$ makes it a unital $C^*$-algebra under the usual $Sup$-norm of functions (Corollary 10).

Then our first major result, Theorem 3, asserts that the $C^*$-algebra $A$ is canonically isomorphic to $A(P(A))$.

With this central theorem (Section 5) in place, the main theme of the present article and its sequels is to deduce results about $C^*$-algebras by considering the corresponding quantum spaces. The following are some such results:

(**1**) In Section 6 we give, as an immediate corollary to Theorem 3, a characterization (Theorem 4) of unitary groups of unital $C^*$-algebras .

(**2**) In Section 7 we give a continuous functional calculus of an *arbitrary* bounded Hilbert space operator. Recall that the Gelfand duality for the unital *commutative* $C^*$-algebra generated by a normal Hilbert space operator $a$ gives a continuous functional calculus for $a$ and can be extended to a Borel functional calculus of $a$. Then, the latter can be given by multiplication operators on a function space, and a special case of this is the Spectral Theorem [27]. When $a$ is not assumed normal, the $C^*$-algebra generated by $\{1, a, a^*\}$ is noncommutative, and our noncommutative Gelfand-Naimark duality leads to a noncommutative continuous functional calculus for $a$, which is realized by noncommutative multiplication operators (Theorem 6).

(**3**) Section 8 is on the Invariant Subspace Theorem, which asserts the existence of nontrivial invaraint subspaces for an arbitrary bounded operator on a several-dimensional complex Hilbert space (Theorem 7) . The theorem has been hitherto proved in the case of normal operators, and for several increasingly general classes of operators (See [19] for a resumé). The result also holds for compact operators [21], and holds trivially for arbitrary operators on *nonseparable* Hilbert spaces. The case of an arbitrary operator on the infinite dimensional *separable* Hilbert space had remained open. Our proof (Theorem 7) essentially covers the separable case—the finite dimensional as well as infinite dimensional.

(**4**) In Section 9, we describe two more applications, the proofs of which will appear in sequels to the present article: (i) Theorem 3 can be easily extended to non-unital $C^*$-algebras. In this case, $P(A)$ is a locally compact Hausdorff q-space. Using this extension, the Pontryagin Duality Theorem for abelian locally compact groups can be extended to arbitrary locally compact groups (Theorem 8). The table in Section 9 gives a quick overview of this result. (ii) An extension of Stone's representation of Boolean algebras to orthomodular lattices. Strictly speaking, this is not an application of the results in the present article. Rather, it is an application of the main ideas surrounding Theorem 3 to an analogous problem in the field of orthomodular lattices.

Finally, several numbered remarks throughout the article point out how various results presented here reduce to standard results in commutative and/or finite-dimensional cases.

## 1.2  Comparision with other approaches

Now, a few words on the works [3, 4, 24] and a comparision with the present work.

(**1**) A theorem of Kadison [17] says that given a $C^*$-algebra $A$, the space of real-valued continuous affine functions on the state space $S(A)$ is isometrically isomorphic to the real space $A_{sa}$ of self adjoint elements of $A$. While this is sufficient to determine the complex Banach space structure of $A$, it does not uniquely determine the multiplicative structure of $A$. Additonal structure, namely, a choice of a continuously varying assignment of orientation to all 3-ball faces of the $S(A)$ uniquely determines the multiplicative structure of $A$ [3].

(**2**) There is another elaboration of this approach due to Shultz [24], which defines on the set $P(A)$ of pure states a uniformity, along with an assignment of an 'orientation' of $P(A)$, and a notion of "transition probability". Then, it is shown that $P(A)$ along with this extra structure uniquely determines the $C^*$-algebra structure of $A$ [24].

Using approaches 1 and 2, the algebra $A$ can be reconstructed as follows. First, as mentioned above, the Banach space structure of $A$, as well as the order structure, is captured using affine continuous functions on $S(A)$. Then, the affine structure gives a functional calculus, using which, one can define a Jordan product on this Banach space. Also, on the same space, there is a naturally defined Lie product compatible with the Jordan product. Recall that a Jordan product has the properties of the anticommutator, and a Lie product is essentially a commutator product. Thus, combining these two products, one gets a product which turns out to be the product of $A$ [3, 4].

As we have seen, our approach reconstructs the algebra $A$ directly and explicitly, as the noncommutative $C^*$-algebra $A(P(A))$ of q-continuous functions on $P(A)$, in perfect analogy with the commutative Gelfand-Naimark. As we shall see, in our setup, a great deal of the intuition of topology carries over to the noncommutative case. This proves to be a very useful heuristic guide as to what might be true in the noncommutative context; not only in the case of $C^*$-algebras, where it leads to an optimal structure theorem for arbitrary operators and a proof of the Invariant Subspace Theorem, but also in other "noncommutative" contexts, such as orthomodular lattices, Pontryagin duality etc.

Another interesting point emerges from the comparision of our approach to 2. The structure in the latter is essential to $C^*$-algerbaic formulation of quantum physics, where states, representations, and transition probabilities play crucial conceptual roles. In confirmity with this conceptual framework, the formulation 2 is probably the most relevant one. However, there are some serious, and virtually insurmountable, conceptual and technical problems associated with such formulations of quantum theory. In view of this, our approach may have some relevance, since it can be viewed as giving a novel (quantum) topological interpretation of the $C^*$-algebra formulation of quantum theory. For example, in this new setup, we can view operators as complex valued functions on a Hausdorff q-space, whence the noncommutativity of the quantum variables appear as a manifestation of the *nonlocal* nature of the product of these functions via the topological structure of the underlying q-space (See Remark 4). This viewpoint may point to a necessary change in the essentially local formulations of quantum field theory which are given in terms of a sort of 'co-sheaf' (a *net*) of $C^*$-algebras [16]. Here, algebras of operators, or more generally, algebras of operator-valued functions, can be replaced by algebras of complex-valued functions which multiply nonlocally—and hence noncommutatively. This seems natural, because, in the physical world, values of these functions are more directly familiar objects, and it is such functions that we would want to be able to multiply without losing the essential noncommutativity of the situation. Furthermore, our viewpoint appears to be *prima facie* evidence that problems of quantum field theory arise because it attempts to fit essentially nonlocal objects into a local strait-jacket.

## 2 PRELIMINARY DISCUSSION

Let $A$ be a unital $C^*$-algebra. A state $\alpha$ on $A$ is pure if and only if the corresponding Gelfand-Naimark-Segal (GNS) representation $\pi_\alpha$ is irreducible [18]. Let $S(A)$ be the set of states on $A$, and $P(A)$ the set of pure states on $A$. Then $S(A)$ is compact Hausdorff in the weak*- topology. However, $P(A)$ is not weak*-closed in general, and hence is not weak*-compact in general. Now we say that $\alpha, \beta \in S(A)$ are **equivalent** if and only if the GNS representations corresponding to $\alpha$ and $\beta$ are equivalent. We denote this equivalence relation by $RS(A)$, and the restriction of this equivalence relation to $P(A)$ will be denoted by $R(A)$.

### 2.1

The following proposition characterizes commutative $C^*$-algebras and is our point of departure.

**Proposition 1** *A unital $C^*$-algebra $A$ is commutative if and only if the equivalence relation $R(A)$ is discrete, i.e. all its equivalence classes are singleton sets. In this case, $R(A) = diag(P(A) \times P(A))$.*

**Proof**: The proof is trivial. At any rate, we are not going to use this result in what follows. Indeed, it is an immediate corollary to Theorem 3. ∎

### 2.2

Let $C(P(A))$ be the algebra of continuous complex-valued functions on $P(A)$. Let

$$A \to C(P(A)) : a \mapsto \widehat{a}$$

be given by

$$\widehat{a}(\alpha) := \alpha(a).$$

We denote the image of this map by $\widehat{A}$. Note that the map $A \to \widehat{A} : a \mapsto \widehat{a}$ is continuous, linear, one-to-one, and preserves the units and the involutions. Also, $\widehat{A}$ is a closed self-adjoint linear subspace of $C(P(A))$ containing $1$, and separates points of $P(A)$.

However, we emphasize that in general, the map $A \to \widehat{A} \hookrightarrow C(P(A))$ is not an algebra homomorphism; the image $\widehat{A}$ is *not* a subalgebra of $C(P(A))$ in general. Indeed, it is a well-known fact [18] that

**Proposition 2** *The following are equivalent:*

1. *$A$ is commutative.*

2. *The map $a \mapsto \widehat{a}$ is a $C^*$-algebra homomorphism.*

3. *The map $a \mapsto \widehat{a}$ is onto, i.e. $\widehat{A} = C(P(A))$.*

4. $\widehat{A}$ is a C*-subalgebra of $C(P(A))$. ∎

Thus, we have identified two obstructions to extending the Gelfand-Duality to noncommutative algberas: in the general non-commutative case, $P(A)$ is not weak*-compact and $\widehat{A}$ is not an algebra. Both these problems can be overcome by considering a natural noncommutative version of topology on $P(A)$ (See Definition 7 and Subsection 5.2). Then, $P(A)$ is a compact Hausdorff in this 'q-topology', and $\widehat{A} = A(P(A))$ where the right-hand side is the set of complex-valued functions on $P(A)$ which are continuous with respect to this topology. Furthermore, this topology gives rise to a natural *noncommutative* product of functions in $A(P(A))$ with respect to which $A(P(A))$ is a $C^*$-algebra isomorphic to $A$.

In next two sections we introduce basic notions of noncommutative topological spaces, leading to the Gelfand duality for arbitrary $C^*$-algebras. We will call these spaces *quantum spaces*.

## 3   QUANTUM SETS

Here we define a notion of quantum sets on which we will add more structure in the next section.

### 3.1   Orthomodular Lattices

**Definition 1** *By an **orthomodular lattice** we shall mean a lattice $L(\wedge, \vee, 0, 1)$ along with a unary operation $\perp\colon p \mapsto p^\perp$, called an **orthocomplementation**, such that the following conditions are satisfied:*

1. *If $p \leqslant q$, then $q^\perp \leqslant p^\perp$.*
2. $(p^\perp)^\perp = p$.
3. $p^\perp \vee p = 1$ *and* $p^\perp \wedge p = 0$.
4. *If $p \leqslant q$, then $p \vee (p^\perp \wedge q) = q$.*

**Remark 1** *We note that the last condition is a weakening of distributive property. An orthomodular lattice satisfying the stronger distributive property is a Boolean algebra. The following proposition shows that in a very precise sense, orthomodular lattices constitute a noncommutative generalization of Boolean algebras.*

**Proposition 3** *Let $L$ be an orthomodular lattice, and $\forall p, q \in L$, define $p \dot\wedge q := p \wedge (p^\perp \vee q)$. Then $L$ is a Boolean algebra if and only if $\forall p, q \in L$, $p \dot\wedge q = q \dot\wedge p$.*

**Proof:** See [5]. ∎

Besides Boolean algebras, the lattice of projections on a Hilbert space, and more generally, lattices of projections of von Neumann algebras, are very important examples of orthomodular lattice (OML). We can define a lattice structure on the set of projections of a von Neumann algebra as follows: for projections $p, q$ we define $p \leq q$ if $p = pq = qp$. This partial order defines an orthomodular lattice structure on projections. In Section 5), another class of OML will be introduced.

## 3.2 Quantum sets and functions

**Definition 2** *Let $X$ be a set and let $P(X)$ be the powerset of $X$. A family $L(X)$ of subsets of a set $X$ which is an orthomodular lattice with respect to the order given by set inclusion is called an **OML of subsets** of $X$. For an OML $L(X)$ of subsets of $X$, the inclusion $L(X) \subset P(X)$ preserves the order. Furthermore, for a family $U_\lambda$ of members of $L(X)$, $\bigwedge_\lambda U_\lambda = \bigcap_\lambda U_\lambda$. An **orthomodular set** is a pair $(X, L(X))$ where $X$ is a set and $L(X)$ is a complete OML of subsets of $X$, such that $X \in L(X)$, and $\emptyset \in L(X)$. By a **quantum set** (or **q-set**) we will mean an orthomodular set $(X, L(X))$, such that for all $U, V, U_\lambda \in L(X)$,*

1. $\emptyset \in L(X)$, and $X \in L(X)$.
2. For all $x \in X$, $\{x\} \in L(X)$.
3. $U \leqslant V$ if and only if $U \subset V$.
4. $\bigwedge_\lambda U_\lambda = \bigcap_\lambda U_\lambda$.
5. $\bigvee_\lambda U_\lambda = \{x \in X : \{x\} \leqslant \bigvee_\lambda \{u_\lambda\}, u_\lambda \in U_\lambda\}$.
6. For all families $\{u_\lambda\} \subset U$, if $\{x\} \leqslant \bigvee_\lambda \{u_\lambda\}$, then $x \in U$.

**Remark 2** *Clearly for any set $X$, $(X, P(X))$ is a quantum set. However, there are three important differences between an arbitrary quantum set $(X, L(X))$ and the **classical set** $(X, P(X))$ :*

1. *We emphasize that $L(X)$ is not distributive in general. Indeed, $L(X)$ is distributive if and only if $L(X)$ is a Boolean algebra if and only if $L(X) = P(X)$.*

2. *Note that $U \cup V \subset U \vee V$, but the reverse inclusion may not hold in general. Indeed, $U \vee V = U \cup V$ if and only if $U \dot\wedge V = V \dot\wedge U$, so that $\vee$ and $\cup$ coincide on $L(X)$ if and only if the latter is a Boolean algebra. In that case, $L(X) = P(X)$.*

3. *Also note that we are not requiring that the orthocomplementation of $L(X)$ co-incide with that of $P(X)$. Indeed, since the wedge operation is simply the set intersection, orthocomplementation of $L(X)$ co-incides with that of $P(X)$ if and only if $\vee$ co-incides with $\cup$ if and only if $L(X) = P(X)$.*

In Section 5 we will meet a large class of quantum sets $(X, L(X))$ where $L(X)$ is not Boolean. These examples come from noncommutative $C^*$-algebras: the set $P(A)$ of pure states of a $C^*$-algebra $A$ constitutes a quantum set, and is a classical set if and only if $A$ is commutative.

**Definition 3** *Given a quantum set $(X, L(X))$ sets $Y \in L(X)$ will be called **quantum subsets** (or **q-subsets**) of $X$.*

It is a standard fact [5] of OML theory that for $Y \in L(X)$ the set $L(Y) := \{U \in L(X) : U \leqslant Y\}$ is an orthomodular lattice with complementation $U \mapsto U^{\perp'}$ defined by $U^{\perp'} := Y \wedge U^\perp$. It follows that for a q-subset $Y$ of $X$, $(Y, L(Y))$ is a quantum set.

**Definition 4** *A **quantum map** (or **q-map**) from a q-set $(X, L(X))$ to a q-set $(Y, L(Y))$ is a usual map $f : X \to Y$ that pulls back quantum subsets of $Y$ to quantum subsets of $X$, i.e. if $U \in L(Y)$, then $f^{-1}(U) \in L(X)$.*

When there is no cause of confusion, we will drop the qualifier 'q' from the terms q-set, q-subset, q-map, etc.

### 3.3 Noncommutative product of functions on quantum sets

The crucial link between quantum spaces and noncommutative $C^*$-algebras arises as follows. Given a quantum set $(X, L(X))$, there is a semigroup (with some extra structure) $S(X) \subset P(X)$, which determines and is determined by the OML $L(X)$. Using the product of $S(X)$, we will later construct noncommutative algebras of complex-valued functions. We describe these matters briefly. (See [13, 7] for details and proofs.)

Given an OML $L$, we write $M(L)$ for the semigroup consisting of all the **monotone maps** $\phi : L \to L$, i.e. maps $\phi$ such that $p \leqslant q$ implies $\phi(p) \leqslant \phi(q)$. For each $p \in L$, the map $\phi_p : L \to L$, given by $\phi_p(q) := p \wedge q$, $\forall q \in L$, is called the **Sasaki projection** corresponding to $p$. Then $\forall p \in L$, $\phi_p \in M(L)$. However, the set of Sasaki projections is not closed under composition. Let $S(L)$ be the sub-semigroup of $M(L)$ generated by the set of Sasaki projections. Then the set of Sasaki projections is closed under composition if and only if it equals $S(L)$, if and only if $L$ is a Boolean algebra. Returning to the general case, for each $\phi \in S(L)$, there exists a unique $\phi^* \in S(L)$ such that $\forall p \in L$, $\phi[\phi^*(p^\perp)^\perp] \leqslant p$, and $\phi^*[\phi(p^\perp)^\perp] \leqslant p$. Then clearly, $(\phi^*)^* = \phi$, and $(\phi\psi)^* = \psi^*\phi^*$, i.e., the map $p \mapsto p^*$ is an **involution** on $S(L)$. Also, for a Sasaki projection $\phi_p$, we have a Sasaki projection $\phi_p^\perp := \phi_{p^\perp}$. Now we extend the map $\phi_p \mapsto \phi_p^\perp$ to $S(L)$ by the identity $(\phi\psi)^\perp = \psi^\perp\phi^\perp$. Then, the set $L(S(L)) := \{\phi \in S(L) : \phi^2 = \phi = \phi^*, (\phi^\perp)^\perp = \phi\}$ of **closed projections** in $S(L)$ is precisely the set of Sasaki projections. In connection with remarks above, $L$ is Boolean if and only if $L(S(L)) = S(L)$. Now for $\phi, \psi \in L(S(L))$, define $\phi \leqslant \psi$ if and only if $\phi = \phi\psi = \psi\phi$. Then it can be shown [13, 7] that $L(S(L))$ is an OML, and $p \mapsto \phi_p$ is an isomorphism of OML's: $L \cong L(S(L))$.

Now we translate the above constructions to the OML $L(X)$ of a quantum set $X$. Let $S(L(X))$ be the semigroup generated by the set of Sasaki projections on $L(X)$. Then, as in the preceding paragraph, there is an involution $\phi \mapsto \phi^*$ on $S(L(X))$. Now for $\phi \in S(L(X))$, let $U_\phi := \phi(X) \in P(X)$. Then, $\phi \mapsto U_\phi$ is a one-to-one map. Let $S(X) \subset P(X)$ be the image of this map. Then by transfer of structure from $S(L(X))$ to $S(X)$, we have a product

$$(U, V) \mapsto U * V$$

defined on $S(X)$, such that $U_\phi * U_\psi = U_{\phi\psi}$, an invoution

$$U \mapsto U^*$$

such that $U_{\phi^*} = (U_\phi)^*$, and $U \mapsto U^\perp$ such that $(U * V)^\perp = V^\perp * U^\perp$. Now the set of closed projections of $S(X)$ equals the set $L(X)$, and defining $U \leqslant V$ if $U = U * V = V * U$ gives an OML structure on $L(X)$ which is the same OML structure that we started with. Thus the OML $L(X)$ completely determines and is determined by the semigroup $S(X)$. We summarize in the following proposition:

**Proposition 4** *Given a q-set $(X, L(X))$, there is a subset $S(X) \subset P(X)$ with an associative product $(U,V) \mapsto U * V$, and an involution $U \mapsto U^*$ and an involutive map (complementation) $U \mapsto U^\perp$ defined on it, making it an involutive semigroup (monoid, actually), such that the lattice of closed projections of $S(X)$ co-incides with the lattice $L(X)$.* ∎

**Remark 3** *Note that $U \wedge V \neq U * V$ in general. Indeed, it can be shown that $U * V = U \wedge V$ if and only if $U * V = V * U$. Thus, $S(X)$ is commutative if and only if $L(X)$ is a Boolean algebra. In this case, $L(X) = S(X)$.*

**Definition 5** *Now we define a **noncommutative product of functions** on a quantum set $(X, L(X))$ as follows. Let $\chi_W$ be the characteristic function of any $W \in S(X)$. For $U, V \in S(X)$ define*

$$\chi_U * \chi_V := \chi_{(U*V)}.$$

*Also, we use the lattice $L(X)$ to define an OML structure on the set of functions of the form $\chi_W, W \in L(X)$ as follows:*

1. $(\chi_U)^\perp := \chi_{U^\perp}$,
2. $\chi_U \wedge \chi_V := \chi_{(U \wedge V)}$,
3. $\chi_U \vee \chi_V := \chi_{(U \vee V)}$,
4. $\mathbf{0} :=$ *The constant function* $0$, *and* $\mathbf{1} :=$ *The constant function* $1$.

*Let $\mathbf{V(X)}$ be the algebra generated by the semigroup of characteristic functions of sets in $S(X)$, and let $\mathbf{F(X)}$ be the uniform closure of the algebra $V(X)$. Then $*$ extends to a product on $F(X)$. We denote this product, too, by $*$.*

**Remark 4** *Since $U * V$ may not be equal to $U \cap V$, the product $*$ of functions is essentially nonlocal, by which we mean that the value of a product $f * g$ at a point $x$ may depend on values of $f$ and $g$ at points other than $x$. Thus the product depends on which quantum set we are taking the product over. If $Y \subset X$, is a quantum subset, and $*_X$ and $*_Y$ the product of functions on $X$ and $Y$ respectively, then $(f|_Y) *_Y (g|_Y) \neq (f *_X g)|_Y$ in general. However, the next proposition asserts that in some cases $(f|_Y) *_Y (g|_Y) = (f *_X g)|_Y$:*

**Definition 6** *A q-subset $Y$ of a q-set $X$ is called **saturated** if for all $U, V \in L(X)$, $(U \wedge Y) *_Y (V \wedge Y) = (U *_X V) \wedge Y$.*

**Proposition 5** *If $Y$ is a saturated q-subset of a q-set $X$, then*

$$(f|_Y) *_Y (g|_Y) = (f *_X g)|_Y.$$

∎

**Remark 5** *In particular, if $L(X) = P(X)$, every subset of $X$ is a saturated subset, so $(f|_Y) *_Y (g|_Y) = (f *_X g)|_Y$, for all $Y \in L(X)$, which is just a manifestation of the fact that for all $x$, $(f * g)(x) = f(x)g(x)$ i.e. that $(f * g)(x)$ depends only on values of $f$ and $g$ at $x$.*

## 4    QUANTUM TOPOLOGY

Now we add more structure to quantum sets.

### 4.1    Quantum spaces and continuous q-functions

**Definition 7** *By a **quantum topology** on a quantum set (or more generally, on an orthomodular set) $(X, L(X))$ we mean a subset $\tau$ of the lattice $L(X)$ satisfying the following conditions.*

1. *$X \in \tau$ and $\emptyset \in \tau$.*

2. *If $U, V \in \tau$ then $U \wedge V \in \tau$.*

3. *If $\{U_\lambda\}_{\lambda \in \Lambda}$ is an arbitrary family of sets in $\tau$, then $(\bigvee_{\lambda \in \Lambda} U_\lambda) \in \tau$.*

*A quantum set (resp. an orhtomodular set) $(X, L(X))$ endowed with a quantum topology $\tau$ will be called a **quantum (topological) space** (resp. **orthomodular space**) and will be denoted by $(X, L(X), \tau)$ or $(X, L(X))$, or even $X$. When $L(X) = P(X)$, $\tau$ is a topology in the usual sense of the word. We shall sometimes refer to such a topology by **classical topology**.*

Now we can extend the entire vocabulary of topology to our more general setup. So we will talk about **quantum open sets, Hausdorff quantum spaces, compact quantum spaces, quantum Borel sets,** etc. We will also use abreviated terms such as **q-set, q-space,** etc. Similarly, we define **q-Borel measures** on a q-space $(X, L(X))$ to be $\sigma$-additive complex-valued functions on the **q-Borel algebra** of $(X, L(X))$. Also, when there is no chance of confusion, we will drop the qualifier 'quantum' or 'q'.

**Definition 8** *A quantum map $f : (X, L(X)) \to (Y, L(Y))$ between quantum spaces will be called **q-continuous** or, briefly, **continuous** if it pulls back q-open sets in $L(Y)$ to q-open sets in $L(X)$.*

**Remark 6** *We emphasize that the notion of a quantum space is a generalization of the classical notion of a topological space. It is this larger category of compact Hausdorff q-spaces that will supply us with 'duals' of unital $C^*$-algebras.*

**Definition 9** *Let $(Y, L(Y))$ be a quantum subset of a q-space $(X, L(X), \tau)$. Then the inclusion $i : Y \hookrightarrow X$ is a quantum map. Let $\tau'$ be the q-topology generated on $(Y, L(Y))$ by the set*

$$\{Y \wedge U : U \in \tau\}$$

*so that $(Y, L(Y), \tau')$ is a quantum space. Then $\tau'$ is the weakest q-topology on $(Y, L(Y))$ making $i$ continuous. The q-set $(Y, L(Y))$ endowed with the topology $\tau'$ will be called a **quantum subspace** of $(X, L(X), \tau)$. When $Y$ is a saturated quantum subset, $\tau'$ is actually equal to $\{Y \wedge U : U \in \tau\}$*

Let $(X, L(X))$ be a q-set and let $\{U_\lambda\}_{\lambda \in \Lambda}$ be a family of elements of $L(X)$. Let $\tau_q$ be the quantum topology generated by $\{U_\lambda\}_{\lambda \in \Lambda}$ and let $\tau_c$ be the classical topology generated by $\{U_\lambda\}_{\lambda \in \Lambda}$. Then, since $\bigcup_i U_i \subset \bigvee_i U_i$ and $U \cap V = U \wedge V$, we see that for every point $x \in X$, every $\tau_q$-neighborhood contains a $\tau_c$-neighborhood, so that $\tau_q \subset \tau_c$. Thus, the classical topology generated by a family of elements of $L(X)$ is in general finer than the quantum topology generated by the same family of subsets. More generally, let $(X, L(X))$ and $(X, L'(X))$ be quantum sets. Then, if $L'(X) \subset L(X)$, and if $\tau$ and $\tau'$ are topologies generated on $(X, L(X))$ and $(X, L'(X))$ respectively by the same family of sets in $L'(X)$, then $\tau'$ is finer than $\tau$.

### 4.2 Product of q-continuous functions

**Proposition 6** *Let $(X, L(X), \tau)$ be a quantum space, and $F(X)$ the algebra of functions defined in Subsection 3.3. Then*

1. *Every $\mathbb{C}$-valued q-continuous function belongs to $F(X)$.*

2. *If $f, g \in F(X)$ are q-continuous, then the product $f * g$ is q-continuous.*

**Proof:**

1. Since q-continuous functions on $X$ can be uniformly approximated by simple functions in $V(X)$, every q-continuous function on $X$ is in $F(X)$.

2. By choosing nets $f_\lambda$ and $g_\lambda$ in $V(X)$ uniformly converging to $f$ and $g$ respectively, $f_\lambda * g_\lambda$ (uniformly) converges to $f * g$. Now considering a net $x_\alpha$ in $X$ converging to $x$, it follows by $\frac{\sqrt{\epsilon}}{3}$ type argument that $f * g$ is continuous. ∎

For a compact Hausdorff quantum space $X$, the set $\boldsymbol{A(X)}$ of q-continuous functions on $X$ is a vector space under pointwise sum and usual scalar product. Then, by Subsection 4.1, $A(X) \subset C(X)$, but as discussed in Section 2, $A(X)$ is not necessarily a subalgebra of $C(X)$. However, the preceding proposition says that $A(X)$ is an algebra with the product $*$. Also, $A(X)$ carries a natural involution given by complex conjugation, and the constant function $1$ is the unit of $A(X)$. Finally, we endow $A(X)$ with the usual $Sup$-norm of functions. Then, we have the following result:

**Proposition 7** *Given a compact Hausdorff q-space $X$, the (noncommutative) algebra $A(X)$ of complex-valued q-continuous functions on $X$ is a unital $C^*$-algebra.*

**Proof:** It remains to show that $A(X)$ is a Banach algebra, and that $\|f * f^*\| = \|f\|^2$. It is clear that the estimate $\|f * g\| \leq \|f\| \|g\|$ holds, and that $A(X)$ is complete under the sup-norm.

Since $\|f * g\| \leq \|f\| \|g\|$, to show $\|f * f^*\| = \|f\|^2$, it suffices to show that $\|f * f^*\| \geq \|f\|^2$. But this is immediate because $|f * f^*| \geq |f f^*|$, which implies that $\|f * f^*\| \geq \|f f^*\| = \|f\|^2$. ∎

**Remark 7** *If $(X, L(X))$ is a classical space, i.e. if $L(X) = P(X)$, then $A(X) = C(X)$, the $C^*$-algebra of continuous complex-valued functions on $X$, with respect to pointwise product.*

**Remark 8** *We remark that certain useful theorems of topology also hold in the context of quantum topology. In these cases, most proofs of the corresponding classical theorems carry over verbatim to the quantum versions. We will mention such theorems as needed.*

## 5 NONCOMMUTATIVE GELFAND-NAIMARK DUALITY

In the following, we will use notation from Section 2. In particular, $S(A)$ is the set of states and $P(A)$ is the set of pure states on a unital $C^*$-algebra $A$.

### 5.1 The q-space $S(A)$ of states

First we define a q-set structure on the set $S(A)$ of states on $A$ as follows.

Given $\alpha, \beta \in S(A)$, let $\langle \alpha, \beta \rangle \subset A^*$ be the subset consisting of elements of the form $c_1 \alpha + c_2 \beta$, $c_1, c_2 \in \mathbb{C}$, and $|c_1|^2 + |c_2|^2 = 1$. Now we define

$$\{\alpha\} \vee \{\beta\} := \begin{cases} \{\alpha, \beta\} & \text{if } (\alpha, \beta) \notin RS(A) \\ S(A) \cap \langle \alpha, \beta \rangle & \text{if } (\alpha, \beta) \in RS(A), \end{cases}$$

and

$$L(S(A)) := \{U \subset S(A) : (\alpha, \beta \in U) \Rightarrow (\forall \text{ states } \gamma \in \{\alpha\} \vee \{\beta\}, \gamma \in U)\}.$$

This is Condition 6 in Definition 2.

Now we define an orthocomplementation on $L(S(A))$. Let $W = A^{**}$ be the enveloping von Neumann algebra of $A$. Then for each projection $p \in W$, there is a net $p_\lambda \in A$ weak*-converging to $p$. We say that a state $\alpha$ is **orthogonal** to a state $\beta$ if there exists a projection $p$ in $W$ such that for a net $p_\lambda \in A$ converging to $p$, $\lim \alpha(p_\lambda) = 1$, and $\lim \beta(p_\lambda) = 0$. Then $q_\lambda := 1 - p_\lambda$ is a net in $A$ converging to the projection $1 - q \in W$ with $\lim \beta(q_\lambda) = 1$, and $\lim \alpha(q_\lambda) = 0$. Thus, orthogonality is a symmetric relation. We denote this relation by $\perp$, i.e., $\alpha \perp \beta$ means that $\alpha$ and $\beta$ are orthogonal. Now define a unary operation $L(S(A)) \to L(S(A)) : U \mapsto U^\perp$ by $U^\perp := \{\alpha \in S(A) : \forall \beta \in U, \beta \perp \alpha\}$.

Now it is easily checked that if $L(S(A))$ is ordered by set inclusion, then $U \mapsto U^\perp$ defined above is an orthocomplementation, and that $L(S(A))$ forms a complete orthomodular lattice of subsets of $P(S(A))$ of $S(A)$, such that for all $U, V, U_\lambda \in L(S(A))$:

1. $\emptyset \in L(S(A))$, and $S(A) \in L(S(A))$.

2. If $\alpha \in S(A)$, then $\{\alpha\} \in L(S(A))$.

3. $U \leqslant V$ if and only if $U \subset V$.

4. $\bigwedge_\lambda U_\lambda = \bigcap_\lambda U_\lambda$.

5. $\bigvee_\lambda U_\lambda = \{\gamma \in S(A) : \gamma \in \vee_\lambda \{\alpha_\lambda\}, \alpha_\lambda \in U_\lambda\}$.

Thus, it follows that $(S(A), L(S(A)))$ is a quantum set. Furthermore, it follows directly from the definition of $\{\alpha\} \vee \{\beta\}$ that if $A$ is commutative, i.e. if $R(A)$ is a discrete equivalence relation, then $L(S(A)) = P(S(A))$ so $(S(A), L(S(A)))$ is just the classical set $S(A)$.

Similarly, we define a q-set structure on $P(A)$ by setting

$$\{\alpha\} \vee \{\beta\} := \begin{cases} \{\alpha, \beta\} & \text{if } (\alpha, \beta) \notin R(A) \\ P(A) \cap \langle \alpha, \beta \rangle & \text{if } (\alpha, \beta) \in R(A) \end{cases}$$

and

$L(P(A)) := \{U \subset P(A) : (\alpha, \beta \in U) \Rightarrow (\forall \text{ pure states } \gamma \in \{\alpha\} \vee \{\beta\}, \gamma \in U)\}$.

Again, $L(P(A))$ forms an orthomodular lattice of subsets of $P(A)$ such that for all $U, V, U_\lambda \in L(P(A))$,

1. $\emptyset \in L(P(A))$, and $P(A) \in L(P(A))$.

2. If $\alpha \in P(A)$, then $\{\alpha\} \in L(P(A))$.

3. $U \leqslant V$ if and only if $U \subset V$.

4. $\bigwedge_\lambda U_\lambda = \bigcap_\lambda U_\lambda$.

5. $\bigvee_\lambda U_\lambda = \{\gamma \in P(A) : \gamma \in \vee_\lambda \{\alpha_\lambda\}, \alpha_\lambda \in U_\lambda\}$.

Also, define $U \mapsto U^\perp$ as in the case of $S(A)$ above. Then $(P(A), L(P(A)))$ is a quantum set. Furthermore, it follows directly from the definition of $\{\alpha\} \vee \{\beta\}$ that when $A$ is commutative, $L(P(A)) = P(P(A))$ so that $(P(A), L(P(A)))$ is just the classical set $P(A)$. It also follows from the definitions of $L(S(A))$ and $L(P(A))$ that $(P(A), L(P(A)))$ is a quantum subset of $(S(A), L(S(A)))$.

We note that each function $\widehat{a} \in \widehat{A}$ is a quantum map from $S(A)$ (resp.P(A)) to $\mathbb{C}$.

Now let $\tau$ be the quantum topology on $(S(A), L(S(A)))$ generated by inverse images of open sets in $\mathbb{C}$ by functions in $\widehat{A}$, i.e., the smallest quantum topology on $(S(A), L(S(A)))$ with respect to which all elements of $\widehat{A}$ are continuous. Let $A(S(A))$ be the set of functions on $S(A)$ continuous with respect to $\tau$.

**Proposition 8** *The set of $S(A)$ endowed with the quantum topology $\tau$ is a compact Hausdorff q-space.*

**Proof:** We know that $S(A)$ is compact in the weak*-topology. The latter being finer than $\tau$ by 4.1, it follows that $(S(A), L(S(A)), \tau)$ is compact. That it is Hausdorff is immediate from the fact that $\widehat{A}$ separates points of $S(A)$. ∎

### 5.2 The q-space $P(A)$ of pure states

As for $S(A)$ above, let $\tau'$ be the quantum topology on $(P(A), L(P(A)))$ generated by inverse images of open sets in $\mathbb{C}$ by functions in $\widehat{A}$, i.e., the smallest quantum topology on $(P(A), L(P(A)))$ with respect to which all elements of $\widehat{A}$ are continuous. Then $\tau'$ coincides with the subspace q-topology inherited from $(S(A), L(S(A)), \tau)$. Let $A(P(A))$ be the set of functions on $P(A)$ continuous with respect to $\tau'$.

**Proposition 9** *A state $\alpha \in S(A)$ is pure if and only if $\forall a, b \in A$*

$$\widehat{(ab)}(\alpha) = (\widehat{a} * \widehat{b})(\alpha).$$

**Proof:**
Let $\alpha$ be a state which satisfies $\forall a, b \in A$, $\widehat{(ab)}(\alpha) = (\widehat{a} * \widehat{b})(\alpha)$. We want to show that $\alpha$ is pure. Let $\beta, \gamma$ be states and let $s, t$ real numbers satisfying $0 < s < 1$, $0 < t < 1$, $s + t = 1$, and $\alpha = s\beta + t\gamma$. Then we will show that $\beta = \gamma$, and hence $\alpha$ is a pure state.

Let $a$ be an hermitian element in $A$. Then the hypothesis, $\widehat{(ab)}(\alpha) = (\widehat{a}*\widehat{b})(\alpha), \forall a,b \in A$, implies that $\widehat{a^2}(\alpha) = (\widehat{a})^2(\alpha)$. Also, the Cauchy-Schwartz inequality for states implies the following:

$$\begin{aligned}
(\alpha(a))^2 &= [\alpha(1.a)]^2 &\leq&\quad \alpha(1^2)\alpha(a^2) = \alpha(a^2) = \widehat{a^2}(\alpha) = (\widehat{a})^2(\alpha), \\
(\beta(a))^2 &= [\beta(1.a)]^2 &\leq&\quad \beta(1^2)\beta(a^2) = \beta(a^2), \\
(\gamma(a))^2 &= [\gamma(1.a)]^2 &\leq&\quad \gamma(1^2)\gamma(a^2) = \gamma(a^2).
\end{aligned}$$

Then it follows that

$$\begin{aligned}
0 &= \widehat{(a^2)}(\alpha) - (\widehat{a})^2(\alpha) \\
&\geq \widehat{(a^2)}(\alpha) - [\alpha(a)]^2 \\
&= \widehat{(a^2)}(s\beta + t\gamma) - [(s\beta + t\gamma)(a)]^2 \\
&= s\beta(a^2) + t\gamma(a^2) - [s\beta(a) + t\gamma(a)]^2 \\
&\geq s[\beta(a)]^2 + t[\gamma(a)]^2 - [s\beta(a) + t\gamma(a)]^2 \\
&= s(s+t)[\beta(a)]^2 + t(s+t)[\gamma(a)]^2 - [s\beta(a) + t\gamma(a)]^2 \\
&= st[\beta(a) - \gamma(a)]^2.
\end{aligned}$$

Thus, $0 \geq st[\beta(a) - \gamma(a)]^2$, and hence, $0 = st[\beta(a) - \gamma(a)]^2$. Consequently, $\beta(a) = \gamma(a)$ for all hermitian $a \in A$. This implies that $\beta = \gamma$ and hence that $\alpha$ is a pure state.

Now we show the converse, *i.e.*, that for each pure $\alpha$, the formula $\widehat{(ab)}(\alpha) = (\widehat{a} * \widehat{b})(\alpha)$ holds. Note that $P(A)$ is a saturated q-subspace of $S(A)$. Hence it follows from Proposition 5 that $\widehat{a} *_{S(A)} \widehat{b} = \widehat{a} *_{P(A)} \widehat{b}$. So it suffices to show the formula with respect to $*_{P(A)}$. To that end, we consider the following. Since every pure state of $A$ can be extended to a pure state of its enveloping von Neumann algebra and since the orthodomodular lattice of projections of a von Neumann algebra completely determines its algebra structure, it suffices to show that for projections $p, q$ in a von Neumann algebra $A$, and for all pure sates $\alpha \in P(A)$,

$$\begin{aligned}
\widehat{p \wedge q}(\alpha) &= (\widehat{p} \wedge \widehat{q})(\alpha) \\
\widehat{p \vee q}(\alpha) &= (\widehat{p} \vee \widehat{q})(\alpha) \\
\widehat{p^\perp}(\alpha) &= (\widehat{p})^\perp.
\end{aligned}$$

Let $U, V \subset P(A)$ such that $\widehat{p} = \chi_U$, and $\widehat{q} = \chi_V$. Then, the above equalities read

$$\begin{aligned}
\widehat{p \wedge q}(\alpha) &= (\chi_U \wedge \chi_V)(\alpha) =: \chi_{U \wedge V}, \\
\widehat{p \vee q}(\alpha) &= (\chi_U \vee \chi_V)(\alpha) =: \chi_{U \vee V}, \\
\widehat{p^\perp}(\alpha) &= (\chi_U)^\perp =: 1 - \chi_U,
\end{aligned}$$

which are easily checked. ∎

**Remark 9** *When $A$ is commutative, $(\widehat{a} * \widehat{b})(\alpha) = \widehat{a}(\alpha)\, \widehat{b}(\alpha)$, by Remark 4. Thus, the equality $\widehat{(ab)}(\alpha) = (\widehat{a} * \widehat{b})(\alpha)$ reads $\widehat{(ab)}(\alpha) = \widehat{a}(\alpha)\, \widehat{b}(\alpha)$, i.e., $\alpha(ab) = \alpha(a)\,\alpha(b)$. Thus, in the commutative case, Proposition 9 simply states the standard fact that a state on $A$ is pure if and only if it is multiplicative.*

**Proposition 10** *The q-space $(P(A), L(P(A)), \tau')$ of pure states is q-compact and q-Hausdorff.*

**Proof:** Let $\alpha_i$ be a net in $P(A)$ converging to $\alpha \in S(A)$. Then, for all $a, b \in A$, $\widehat{(ab)}(\alpha_i) = (\widehat{a} * \widehat{b})(\alpha_i)$. Thus, $\widehat{(ab)}(\alpha) = \widehat{(ab)}(\lim \alpha_i) = \lim \widehat{(ab)}(\alpha_i) = \lim(\widehat{a} * \widehat{b})(\alpha_i) = (\widehat{a} * \widehat{b})(\alpha)$. It follows that $\alpha \in P(A)$. Thus, $P(A)$ is q-closed in $S(A)$, and hence is q-compact. Also, since $\widehat{A}$ separates points of $P(A)$, it follows that $P(A)$ is q-Hausdorff. ∎

**Corollary 1** *For a $C^*$-algebra $A$, the set $A(P(A))$ of q-continuous functions on $P(A)$ is a unital $C^*$-algebra.*

**Proof:** Immediate from Proposition 7 and Proposition 10. ∎

### 5.3 Noncommutative Gelfand-Naimark Duality

We will use the following extension of classical Stone-Weierstrass theorem in the proof the theorem on Noncommutative Gelfand-Naimark Duality.

**Theorem 2 (Noncommutative Stone-Weierstrass Theorem)** *Let $X$ be a q-compact and q-Hausdorff q-space, and let $A(X)$ be the unital $C^*$-algebra of q-continuous complex-valued functions on $X$. Let $B$ be a unital $C^*$-subalgebra of $A(X)$ which separates points of $X$. Then we can conclude that $B = A(X)$.*

**Proof:** The proof is verbatim the same as in the commutative case. See [23], for example. ∎

**Theorem 3 (Noncommutative Gelfand-Naimark Duality)** *For a unital $C^*$-algebra $A$, the map $A \to A(P(A)) : a \mapsto \widehat{a}$ is an isomorphism of unital $C^*$-algebras:*

$$A \cong A(P(A)).$$

**Proof:**
It follows from Proposition 9 that the map is a $C^*$-homomorphism. Let $\widehat{A} \subset A(P(A))$ be the image of the map.

The map is one-to-one, for if $\widehat{a} = 0$, then, for all $\alpha \in P(A)$, we have $\widehat{a}(\alpha) = 0$, i.e., $\forall \alpha \in P(A), \alpha(a) = 0$, and hence $a = 0$.

Now we show that the map is onto. Given $\alpha, \beta \in P(A)$, such that $\alpha \neq \beta$, there is an $a \in A$ such that $\alpha(a) \neq \beta(a)$, i.e., $\widehat{a}(\alpha) \neq \widehat{a}(\beta)$. Thus, $\widehat{A}$ separates points of $P(A)$. Clearly, $\widehat{A}$ is self-adjoint and contains $\widehat{1}$, the unit of $A(P(A))$. Consequently, by Theorem 2, the non-commutative Stone-Weierstrass theorem, $\widehat{A}$ is the whole of $A(P(A))$. ∎

**Remark 10** *When $A$ is commutative, $R(A)$ is the discrete equivalence relation. Consequently, $(P(A), L(P(A)))$ is the classical compact Hausdorff space $P(A)$ with weak$^*$-topology, and $A = A(P(A)) = C(P(A))$. Thus, in the commutative case, we recover the Gelfand-Naimark duality (Theorem 1).*

As an easy corollary to the preceding theorem, we can recover Dauns-Hofmann [9] representation of a $C^*$-algebra $A$ as continuous sections of a 'sheaf' (a field) of (presumably simpler) $C^*$-algebras over the spectrum $Sp(A)$ of $A$. This is achieved via the natural quotient map $P(A) \to Sp(A)$.

## 6 CHARACTERIZING UNITARY GROUPS OF UNITAL $C^*$-ALGEBRAS

Let $U(1) := \{\lambda \in \mathbb{C} : |\lambda| = 1\}$. Then the following characterization of groups of unitary elements of unital $C^*$-algebras is immediate from Theorem 3.

**Theorem 4 (Characterization of unitary groups of unital $C^*$-algberas)** *A topological group is homeomorphically isomorphic to the unitary group of a unital $C^*$-algebra if and only if it is homeomorphically isomorphic to the group of q-continuous $U(1)$-valued functions on a compact Hausdorff q-space.* ∎

## 7 CONTINUOUS FUNCTIONAL CALCULUS FOR HILBERT SPACE OPERATORS

Let $B(H)$ be the algebra of bounded linear operators on a Hilbert space $H$. Let $a \in B(H)$, and let $A$ be the unital $C^*$-subalgebra of $B(H)$ generated by $\{1, a, a^*\}$. When $a$ is normal, $A$ is commutative, and $R(A)$ is discrete. In this case, $P(A) \cong \sigma(a)$, the spectrum of $a$, and the commutative Gelfand-Naimark theorem leads to a simultaneous representation of operators in $A$ in terms of multiplication operators on a Hilbert space $L^2(\sigma(a), \mu)$.

In the general case, where $a$ is not assumed normal, the isomorphism of noncommutative $C^*$-algebras, $A \xrightarrow{\cong} A(P(A))$ from Theorem 3, is represented by noncommutative multiplication operators.

## 7.1 Noncommutative Riesz Representation Theorem

We will need the following noncommutative generalization of the Riesz Representation Theorem.

**Theorem 5 (Noncommutative Riesz Representation Theorem)** *Let $X$ be a compact Hausdorff q-space. Then every bounded linear functional $\alpha$ on $A(X)$ is given by integration with respect to a complex-valued q-Borel measure $\mu_\alpha$ on $X$. Furthermore, the map $\alpha \mapsto \mu_\alpha$ is an isomorphism of Banach spaces*

$$A(X)^* \cong M(X),$$

*where $M(X)$ is the space of complex-valued regular bounded q-Borel measures on $X$.*

**Proof:** The proof carries over verbatim from the commutative case. See, for example, [22] . ∎

**Remark 11** *If $X$ is a classical compact Hausdorff space, then $A(X) = C(X)$, and the preceding theorem reduces to the classical Riesz Representation Theorem.*

## 7.2 Functional Calculus of an arbitrary operator

Now we come to the main theorem of this section. Let $\Sigma(a)$ be the image of the map $\widehat{a} : P(A) \to \mathbb{C}$.

**Theorem 6 (The Spectral Theorem and the Functional Calculus)** *Let $a$ be a bounded linear operator on a Hilbert space $H$, and let $A$ be the $C^*$-algebra of operators on $H$ generated by $\{a, a^*, 1\}$. Then, there exists a q-Borel measure $\mu$ on $P(A)$, q-space $X$, which is a disjoint union of a q-Borel measurable subsets $X_i \subset P(A)$, and a unitary isomorphism $U : L^2(X, \mu) \to H$, such that*

1. *If $a$ has a star-cyclic vector in $H$, i.e., if $H$ has an $A$-cyclic vector, then $X = P(A)$.*

2. *For each $x \in A$, we have $x = U M_{\widehat{x}} U^{-1}$, where $M_{\widehat{x}}$ is the multiplication operator $L^2(X, \mu) \to L^2(X, \mu) : f \mapsto \widehat{x} * f$, i.e., $x$ is unitarily equivalent (via $U$) to the noncommutative multiplication operator $M_{\widehat{x}}$.*

3. *The continuous functional calculus of $a$, i.e., the isomorphism $A(X) \to A : \phi \mapsto \phi(a)$ from Theorem 3, is given by*

$$\phi(a) = U M_\phi U^{-1},$$

*where $M_\phi$ is the multiplication operator $L^2(X, \mu) \to L^2(X, \mu) : f \mapsto \phi * f$, i.e., $\phi(a)$ is unitarily equivalent (via $U$) to the noncommutative multiplication operator $M_\phi$. A special case of the formula is:*

$$a \ = U M_{\widehat{a}} U^{-1}.$$

4. *The functional calculus of $a$, $A(X) \to B(H)$ given above uniquely extends the Holomorphic Functional Calculus (when holomorphic functions on $\Sigma(a)$ are duly pulled back to $X$.)*

**Proof:**
The proof is quite similar to that for Normal Operators [27, 10] with the pointwise multiplication of functions replaced by the noncommutative product $*$.

First note that any bounded operator $a \in B(H)$ can be decomposed as a direct sum $a = \bigoplus_{i \in I} a_i$ where each $a_i$ is a star-cyclic operator on closed subspace $H_i$ of $H$, such that $H = \bigoplus_{i \in I} H_i$ with $\{H_i\}_{i \in I}$ mutually orthogonal. Thus, it suffices to show the theorem for an operator $a$ which is star-cyclic on a Hilbert space $H$. Then we can put together the description $a$ from descriptions of $a_i$. (See [27, 10])

So let $h \in H$ be a star-cyclic vector for $a$. (Note that in this case $H$ is necessarily separable.) Then define a bounded functional $A(P(A)) \mapsto \mathbb{C}$ given by $\phi \mapsto \langle \phi(a)h, h \rangle$. By the noncommutative Riesz Representation Theorem (Theorema 5), this defines a q-measure $\mu$ on $P(A)$. Then it is straightforward to see that $U : A(P(A)) \to H$ given by

$$U(\phi) := \phi(a)h$$

extends to an isomorphism $U : L^2(P(A), \mu) \to H$, with $\|U\phi\| = \|\phi\|$. Now for $\phi, \psi \in A(P(A))$, we have $(U^{-1}\phi(a)U)\psi = U^{-1}\phi(a)(U\psi) = (U^{-1}\phi(a)\psi(a))h = U^{-1}\{(\phi * \psi)(a)\}(h) = \phi * \psi = M_\phi \psi$. Consequently, by density of $A(P(A))$ in $L^2(P(A), \mu)$, we have $(U^{-1}\phi(a)U)\psi = M_\phi \psi$ for all $\psi \in L^2(P(A), \mu)$. ∎

## 8 THE INVARIANT SUBSPACE THEOREM

We are now ready for the Invariant Subspace Theorem.

**Theorem 7 (Invariant Subspace Theorem)** *Every bounded operator on a complex Hilbert space $H$ with $\dim(H) > 1$ has a nontrivial invariant subspace.*

**Proof:** Let $a \in B(H)$. Note that it is sufficient to consider $a$ which is star-cyclic on a Hilbert space $H$, with $\sigma(a)$ being singleton. Recall that $\Sigma(a)$ is the image of the map $\widehat{a} : P(A) \to \mathbb{C}$. Then, it is easy to see that $\sigma(a) \subset \Sigma(a)$. We consider the following two cases.

Case(i) Assume $\sigma(a) = \Sigma(a) = \{\lambda\}$.
Then, for all pure states $\alpha \in P(A)$, we have $\alpha(a) = \lambda$, and hence $\alpha(a - \lambda) = 0$, for all $\alpha \in P(A)$. This is equivalent to $a - \lambda = 0$, i.e., $a = \lambda$ a scalar operator, which always has a nontrvial invariant subspace if $\dim(H) \geq 1$.

Case(ii) Assume $\sigma(a) \neq \Sigma(a)$.
Now, by Part (i) of Theorem 6, $H \cong L^2(P(A), \mu)$ and $a$ is (equivalent to) the noncommutative multiplication operator $M_{\widehat{a}}$. Let $\sigma = \widehat{a}^{-1}(\sigma(a))$, $\sigma'(a) := \Sigma(a) \setminus \sigma(a)$, and $\sigma' = \widehat{a}^{-1}(\sigma'(a))$, Thus $P(A) = \sigma \cup \sigma'$.

Then, $L^2(P(A), \mu) = L_\sigma \oplus L_{\sigma'}$, where $L_\sigma := L^2(\sigma, \mu|_\sigma)$ and $L_{\sigma'} := L^2(\sigma', \mu|_{\sigma'})$. Now, since $\sigma'$ contains a non-empty open set $\sigma^\perp$, we have $L_{\sigma'} \neq 0$. It follows that $L_\sigma \neq L^2(P(A), \mu)$. Also, $L_\sigma \neq 0$ because otherwise $a$ will have empty spectrum, which is not possible. Thus, $L_\sigma$ is a nontrivial subspace, which is invariant for $a$, because $a(L_\sigma) \subset L_\sigma$. ∎

## 9 FURTHER APPLICATIONS

We describe here two more applications. Details will appear elsewhere.

### 9.1 Nonabelian Pontryagin Duality

Recall that the set $\widehat{G}$ of characters of a locally compact abelian group $G$ forms a locally compact abelian group and the celebrated Pontryagin duality theorem gives a natural isomorphism $G \cong \widehat{\widehat{G}}$. Extending this theorem to nonabelian groups necessarilly leads to a new notion:

**Definition 10** *A **quantum group space** is a quantum space with a group structure compatible with the quantum toplogy. In this setting, the terms* 'abelian' *and* 'nonabelian' *will refer to the group structure of a quantum group space, and* 'commutative' *and* 'noncommutative' *will refer to its topology.*

Given a locally compact group $G$ we define (see below) its dual to be a certain quantum group space $\widehat{G}$, which is a group if and only if $G$ is abelian. The classical dual of a possibly nonabelian $G$, i.e. the set of equivalence classes of irreducible unitary representations of $G$, is a quotient space of $\widehat{G}$. In the abelian case, $\widehat{G}$ coincides with the classical dual. This viewpoint inevitably leads to an extension of the duality to *quantum group spaces*.

Let $G$ be a locally compact (Hausdorff) quantum group space. Let $A(G)$ be the $C^*$-algebra of q-continuous complex-valued functions on $G$ vanishing at infinity. Then $A(G)$ has a co-product arising from the group structure of $G$, and its enveloping von Neumann algebra $K^*(G) := A(G)^{**}$ is a von Neumann bi-algebra. Now we can construct from the dual von Neumann bialgebra $\widehat{K^*(G)}$ a locally compact quantum space $\widehat{G}$, which has a multiplication structure derived from the co-multiplication of $\widehat{K^*(G)}$. This makes $\widehat{G}$ a quantum group space which we call **the dual quantum group space of $G$**. Repeating this procedure, we can construct the locally compact quantum group space $\widehat{\widehat{G}}$, i.e., the dual of $\widehat{G}$, from $\widehat{K^*(\widehat{G})}$. Then, the following generalization of classical Pontryagin duality holds:

**Theorem 8 (Pontryagin duality for quantum group spaces)** *For a locally compact quantum group space $G$, the dual $\widehat{G}$ is a locally comapct quantum group space, and*

$$G \cong \widehat{\widehat{G}}.$$

Now, let $G, H, K, N$ be quantum group spaces with the corresponding duals $\widehat{G}, \widehat{H}, \widehat{K}, \widehat{N}$. Then the following table summarizes the various situations covered by Theorem 8:

| ↓ Group-Space → | Commutative | Noncommutative |
|---|---|---|
| Abelian | $G$, $\widehat{G}$ | $\widehat{H}$, $K$ |
| Nonabelian | $H$, $\widehat{K}$ | $N$, $\widehat{N}$ |

Thus, the dual $\widehat{G}$ of an abelian group $G$ is an abelian group; for a nonabelian group $H$, the dual $\widehat{H}$ is an abelian noncommutative quantum group space, etc. We note that the box containing $G, \widehat{G}$ is the classical Pontryagin duality. The boxes containing $H, \widehat{K}$ and $K, \widehat{H}$ include nonabelian groups and abelian noncommutative group spaces, and finally the box containing $N, \widehat{N}$ cover nonabelian noncommutative quantum group spaces.

### 9.2 Stone Duality for Noncommutative Boolean algebras, *i.e.*, Orthomodular Lattices

The ideas of Section 5 can be applied to Orthomodular Lattices (OML) (Definition 1). Elements of a Boolean algebra $B$ are represented by clopen subsets of a totally disconnected compact space—the maximal ideal space of $B$ (Stone's Theorem [25]). As in the case of $C^*$-algebras, the geometric object corresponding to a (possibly non-Boolean) OML is a totally disconnected compact orthomodular space (Definition 7) naturally associated with the lattice. Furthermore, an OML is Boolean if and only if this noncommutative space is a usual topological space. In this case, one recovers Stone's theorem. The general case yields an OML analog of Dauns-Hoffman theorem—the Graves-Selesnick representation [15] of an OML as sections of a sheaf of (presumably simpler) OML's.